\documentclass[12pt,a4paper,english]{article}
\usepackage[utf8]{inputenc}
\usepackage{babel}

\usepackage{graphicx} 
\graphicspath{{./img/}} 
\usepackage{booktabs}
\usepackage[colorlinks]{hyperref}
\usepackage{nicematrix} 
\usepackage{pgfgantt}
\usepackage{optidef} 
\usepackage[noabbrev]{cleveref} 
\usepackage{mathtools}

\usepackage{amsfonts}
\usepackage{amsthm}
\usepackage{siunitx}
\usepackage[kerning,spacing,babel]{microtype}
\usepackage[final]{changes}
\theoremstyle{remark} 
\newtheorem*{ex}{Example}
\newtheorem*{im}{Implementation}
\newtheorem*{ju}{Justification}

\usepackage{enumitem}
\setlist[description]{itemsep=0mm}
\setlist[enumerate]{itemsep=0mm}
\setlist[itemize]{itemsep=0mm}	

\newcommand{\e}[1]{\begin{ex}#1\end{ex}} 
\newcommand{\imp}[1]{\begin{im}#1\end{im}} 
\newcommand{\just}[1]{\begin{ju}#1\end{ju}} 
\newcommand{\q}[1]{``#1''} 

\definecolor{orange}{HTML}{FDA448}
\definecolor{blue}{HTML}{469BA9}
\definecolor{orange}{HTML}{fb8072}
\definecolor{blue}{HTML}{80b1d3}
\definecolor{linkcolor}{HTML}{CA0020}
\definecolor{linkcolor}{rgb}{0,0.2,0.6}
\definecolor{linkcolor}{HTML}{145680}
\definecolor{linkcolor}{HTML}{990004}

\hypersetup{
 linkcolor=linkcolor
,citecolor=linkcolor
,filecolor=linkcolor
,urlcolor=linkcolor
,menucolor=linkcolor
,runcolor=linkcolor
,linkbordercolor=linkcolor
,citebordercolor=linkcolor
,filebordercolor=linkcolor
,urlbordercolor=linkcolor
,menubordercolor=linkcolor
,runbordercolor=linkcolor
}

\title{Fast and exact audit scheduling optimization}
\author{Corresponding author: Jan Motl\\
	\textsc{Czech Technical University in Prague}\\
	Thákurova 9, 160 00 Prague, Czechia\\
	+420 603 885 753, jan.motl@fit.cvut.cz\\
	\\
	Pavel Kordík\\
	\textsc{Czech Technical University in Prague}\\
	Thákurova 9, 160 00 Prague, Czechia\\
	+420 604 499 078, pavel.kordik@fit.cvut.cz}
 
\begin{document}
\date{October 24, 2020}
\maketitle

\begin{abstract}

This article is concerned with the cost and time-effective scheduling of financial auditors with integer linear programming. 
The schedule optimization considers 13 different constraints, staff scarcity, frequent alterations of the input data with the need to minimize the changes in the generated schedule, and scaling issues.  
The delivered implementation reduced the time to the first schedule from 3 man-days to 1 h and the schedule update time from 1 man-day to 4 min.

\end{abstract}


Scheduling is an old optimization problem \cite{Dantzig1954,Johnson1954,Smith1956}. However, the difficulty with the optimization is that we have to carefully select what to optimize and which constraints to consider to get a solution for the optimization problem in an acceptable timeframe. Hence, the formulation of the schedule optimization necessarily differs from \deleted{a} firm to firm as their needs differ. In this article, we list constraints encountered when optimizing a schedule for an auditing firm and we illustrate how to formulate them as a multi-commodity network flow problem \cite{ElAdoly2018}.

An unnamed auditing firm suffers from staff scarcity as the number of clients continues to increase. This has led to severe difficulties in planning activities. The last handmade schedule took 3 man-days to create and each schedule update took another 1 man-day. It is predicted that the firm will continue to grow, and the situation will become even worse. Hence, the  management decided that the firm needs a computer-aided planner, as presented in this article.

Our audit-staff scheduling problem may succinctly be characterized by the following statements:
\begin{enumerate}
  \item The firm employs financial \textit{auditors} for audit \textit{engagements} (an audit engagement is an arrangement that the firm has with a client to conduct an audit of the client) during the following year.
  \item Each engagement includes an availability calendar indicating when the audit can be performed. 
  \item Each engagement consists of \textit{tasks}. Each task takes a specified amount of time. Each task is completed by a single auditor. 
  \item Each auditor has an availability calendar with the number of hours that the auditor can work on that day.
  \item We want to find an assignment of the auditors to the engagements that respects the business rules further listed in \cref{sec:method}.
\end{enumerate}

\begin{figure}[htb]
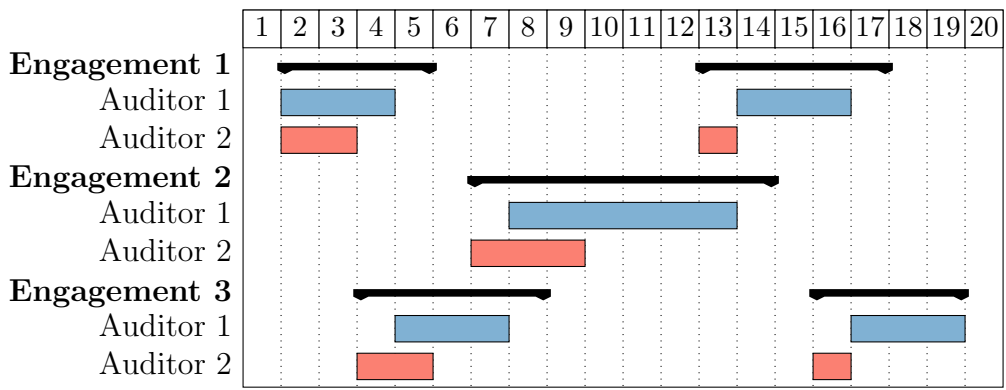

\centering
\begin{ganttchart}[vgrid,
    bar/.append style={fill=gray!50},
    y unit title = 0.5cm, title height=1,
    y unit chart = 0.5cm,
    bar top shift = 0.1, bar height = 0.7
]{1}{20}
\gantttitlelist{1,...,20}{1} \\
   \ganttgroup{Engagement 1}{2}{5}  \ganttgroup{}{13}{17} \\
   \ganttbar[bar/.append style={fill=blue}]{Auditor 1}{2}{4}  
   		\ganttbar[bar/.append style={fill=blue}]{}{14}{16} \\
   \ganttbar[bar/.append style={fill=orange}]{Auditor 2}{2}{3}  
   		\ganttbar[bar/.append style={fill=orange}]{}{13}{13} \\

   \ganttgroup{Engagement 2}{7}{14} \\
   \ganttbar[bar/.append style={fill=blue}]{Auditor 1}{8}{13} \\
   \ganttbar[bar/.append style={fill=orange}]{Auditor 2}{7}{9} \\

   \ganttgroup{Engagement 3}{4}{8}  \ganttgroup{}{16}{19} \\
   \ganttbar[bar/.append style={fill=blue}]{Auditor 1}{5}{7}
   		\ganttbar[bar/.append style={fill=blue}]{}{17}{19} \\
   \ganttbar[bar/.append style={fill=orange}]{Auditor 2}{4}{5} 
   		\ganttbar[bar/.append style={fill=orange}]{}{16}{16}
\end{ganttchart}
\caption{An example of a schedule for two auditors and three engagements. Engagements can be divided into phases, for example, interim and final in engagements 1 and 3. The engagement availability window is depicted within black ranges.}
\label{fig:example}
\end{figure}

A basic illustration of the desired schedule is provided in \cref{fig:example}.

When we formulated the scheduling optimization problem as a multi-commodity network flow problem with 13 constraints, we reduced the time to the first schedule from 3 man-days with a handmade approach to 1 h and the schedule update time from 1 man-day to 4 min.

\paragraph{Paper structure}
First, we \hyperref[sec:review]{review} the audit scheduling literature. We then \hyperref[sec:method]{justify} our choice of the optimization solver because it affects which constraints can be efficiently modeled and how they should be formulated. In the \hyperref[sec:requirements]{requirements} subsection, we list the modeled constraints and provide examples to illustrate their relevance to audit scheduling. In the \hyperref[sec:implementation]{implementation} section, we describe the data model. The \hyperref[sec:results]{result} section then provides an empirical comparison of the optimized schedule versus a handmade schedule, which was created by an experienced scheduler over a 1-month period. Finally, we end the article with a \hyperref[sec:discussion]{discussion} section and concluding remarks in the \hyperref[sec:conclusion]{conclusion} section.

\section{Review} \label{sec:review}
This section reviews different methods that can be used to solve the audit scheduling methods.

Summers \cite{Summers1974} used linear programming (LP) to minimize the cost of an auditor-project assignment based on the level and hourly cost of each auditor and the project's hourly requirements at each level.

Balachandran and Zoltners \cite{Zoltners1989} extended Summer's LP problem formulation to an integer linear programming (ILP) problem formulation to avoid the expressivity limitations of LP.

Romeijn and Morales \cite{Romeijn2000} introduced a class of greedy approximation algorithms for the assignment problem. As an advantage of greedy algorithms, they scale to large instances. However, they do not guarantee the return of a feasible, let alone optimal, schedule. 

Wang and Kong \cite{Wang2012b} used a genetic algorithm (GA) to maximize the match between auditors' specialization and the project characteristics. They chose to use a GA to obtain an approximately optimal solution because they were concerned that obtaining an optimal solution would be prohibitively time-consuming \cite{Wang2012b}[section 2.2].

An extensive review of staff scheduling is provided by Ernst et al. \cite{Ernst2004}, Van Den Bergh et al. \cite{VanDenBergh2013}, and Kalra and Singh \cite{Kalra2015}.

\section{Method} \label{sec:method}
We tested multiple optimization methods, including greedy and genetic algorithms. However, we were unable to obtain a feasible solution from either of these algorithms in real-world instances. We identified two reasons for this. 
\begin{itemize}
  \item Staff scarcity. To cover all the engagements with the limited number of the auditors, the generated schedule must be nearly perfect. For this reason, a pure greedy algorithm is unlikely to provide a valid schedule. This is in agreement with the observation by Drexl and Gruenewald \cite{Drexl1993}.
  \item There is a high count of diverse hard constraints that must be satisfied to obtain a solution acceptable by the auditors. This makes it difficult to design mutation and crossover operators in a genetic algorithm that would maintain the validity of the solution.
\end{itemize}
For these reasons, we opted to use \replaced{an}{the} ILP solver because it is guaranteed to (eventually) find a valid solution if such a solution exists. Once a valid solution is found, all improvements are guaranteed to be valid as well.

\subsection{Formulation}
We denote the sets by capital letter, for example, $A$, and their corresponding variables by lowercase characters, for example, $a$ \added{for} $a \in A$.

We assume that the firm has $A$ auditors, and we want to find a schedule for the auditors for the next $D$ days. The firm has $E$ engagements, where each engagement is divided into consecutive $P$ phases (e.g., interim and final). Each engagement requires auditors at a specific level (e.g., analyst, senior manager). We assume that each auditor is exactly at one level out of all $L$ levels. Because an engagement may require multiple auditors at the same level, we differentiate between them using $I$ indexes.

By introducing aliases, task $T = (E, P, L, I)$ and sink $S = (E, P, L, I, D)$, where day $d$ indicates when task $t$ starts (\cref{eq:searchspace}), we can formulate the auditor schedule optimization as a multi-commodity network flow problem (a multi-commodity network is a network with one or more sources and one or more sinks \cite{ElAdoly2018,Salimifard2020}).

\begin{figure*}[htb]
\begin{equation*}
\text{Auditor} \times \underbrace{\overbrace{\text{Engagement} \times \text{Phase} \times \text{Level} \times \text{Index}}^{\text{Task}} \times \text{Day}}_{\text{Sink}}	
\end{equation*}
	\caption{Search space}
	\label{eq:searchspace}
\end{figure*}

\begin{figure*}[htb!]
	\centering
	\includegraphics[width=0.5\textwidth]{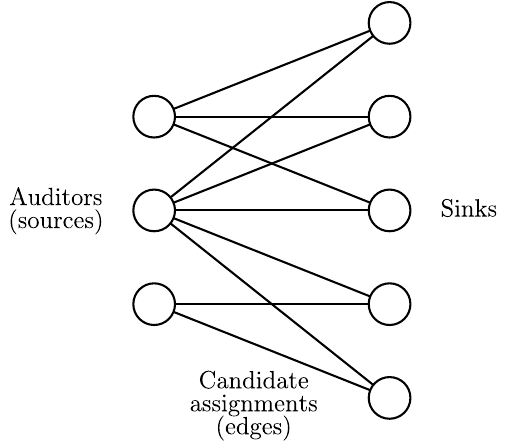}
	\caption{Scheduling problem formulated as a flow problem.}
	\label{fig:flow}
\end{figure*}

The flow problem can be represented as a graph with auditors $A$ on the left and sinks $S$ on the right (see \cref{fig:flow}). 

If auditor $a$ can be potentially assigned to a sink $s$, there is an edge $x_{as}$ between the nodes. The task of the optimizer is then to decide which of the candidate assignments are selected and which are not:
 
\begin{equation}
	x_{as} =
	\begin{cases}
		1 & \text{if auditor $a$ is assigned to sink $s$,} \\
	    0 & \text{otherwise,}
	\end{cases}
\end{equation}
\replaced{where}{For the reminder,} the sink $s$ tells us which task the auditor $a$ should do and on which day $d$ auditor $a$ should start working.

The objective of the optimizer is to minimize the sum of the costs $c_{as}$ associated with assigning auditor $a$ to sink $s$:

\begin{mini}
{x}{\sum_{a \in A} \sum_{s \in S} c_{as} x_{as}}
{\label{eq:Example1}}{}
\addConstraint{\sum_{a \in A} \sum_{d \in D} x_{as}}{ = 1 \text{ for each task $t$},}
\end{mini}
where the constraint states that each task must be assigned. The following further describe\added{s} the basic assignment model. 

\subsection{Requirements}\label{sec:requirements}
Considerable effort was made to extract a complete set of optimization criteria and constraints from the auditors. We provide an exhaustive list of the auditor requirements in the hope that interested readers may take inspiration from them for their own needs. To facilitate this, we do not only list the requirements, but we also justify the requirements, provide illustrative examples and implementation details.

\paragraph{Multitasking}
Based on the poor previous experience with auditors working on multiple tasks in a single day, it was required that a single auditor be assigned at most a single task per day. 
\imp{
We enforced this with the following constraint:
\begin{equation}
	\sum_{t \in O_{ad}} x_{as}\leq 1 \text{ for each auditor $a$ and day $d$},
\end{equation}
}
\noindent
where $O_{ad}$ is a set of overlapping candidate tasks for auditor $a$ and day $d$.

\paragraph{Continuity}\phantomsection \label{p:continuity}
It was also required that an auditor finish\added{es} one task before moving to another to minimize the context switch penalty.
\imp{
We solve this by generating only edges, where an auditor utilizes all its available time to the potential task until the task is finished (see \cref{fig:edges}). \added{The task durations are given by the engagement managers and the durations are assumed to be the exact estimates, which are independent of the assigned auditor.}

\begin{figure}[htb]
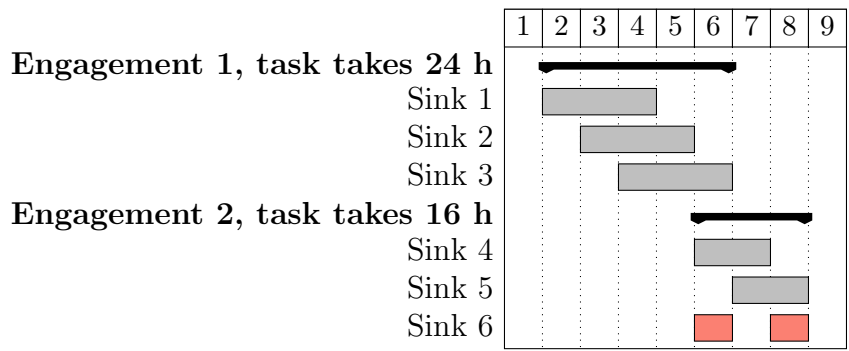

\centering
\begin{ganttchart}[vgrid,
    bar/.append style={fill=gray!50},
    y unit title = 0.5cm, title height=1,
    y unit chart = 0.5cm,
    bar top shift = 0.1, bar height = 0.7
]{1}{9}
  \gantttitlelist{1,...,9}{1} \\
  \ganttgroup{Engagement 1, task takes 24 h}{2}{6} \\
  \ganttbar{Sink 1}{2}{4} \\
  \ganttbar{Sink 2}{3}{5} \\
  \ganttbar{Sink 3}{4}{6} \\ 
  
  \ganttgroup{Engagement 2, task takes 16 h}{6}{8} \\
  \ganttbar{Sink 4}{6}{7} \\
  \ganttbar{Sink 5}{7}{8} \\ 
  \ganttbar[bar/.append style={fill=orange}]{Sink 6}{6}{6} 
  \ganttbar[bar/.append style={fill=orange}]{}{8}{8} 
\end{ganttchart}
\caption{Availability of each engagement in days is depicted by the black ranges. Possible sinks, assuming that auditor $a$ can work 8 h every day, are represented by gray rectangles. Sinks 3 and 4 cannot be assigned to auditor $a$ because of the overlap. However, sinks 2 and 4 can be assigned to auditor $a$ because no overlap occurs. An example of an illegal sink is sink 6, because once an auditor starts working on an engagement, the auditor must finish it as soon as possible.}
\label{fig:edges}
\end{figure}
}

\paragraph{Auditor availability}\phantomsection \label{p:auditorAvailability}
Each auditor has a calendar with the number of hours they can work on a specified day. This allowed us to model the following constraints:
\e{The auditor will be hired at the beginning of the next month.}
\e{It is known that the auditor will leave the firm in 2 months.}
\e{The auditor works part-time.}
\e{Training, vacations, and national holidays.}
\imp{
When the auditor does not work on day $d$, the auditor cannot start a new task on that day. This is ensured by not generating edges in the graph.

We also skip the generation of edges for tasks that auditor $a$ may not finish in time (see \cref{fig:availability}).

\begin{figure}[htb]
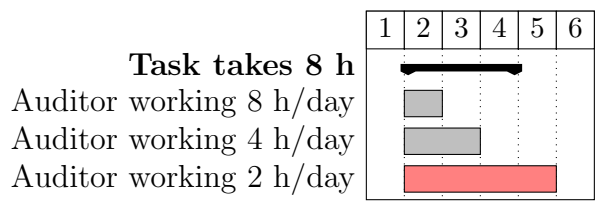

\centering
\begin{ganttchart}[vgrid,
    bar/.append style={fill=gray!50},
    y unit title = 0.5cm, title height=1,
    y unit chart = 0.5cm,
    bar top shift = 0.1, bar height = 0.7
]{1}{6}
  \gantttitlelist{1,...,6}{1} \\
  \ganttgroup{Task takes 8 h}{2}{4} \\
  \ganttbar{Auditor working 8 h/day}{2}{2} \\
  \ganttbar{Auditor working 4 h/day}{2}{3} \\
  \ganttbar[bar/.append style={fill=red!50}]{Auditor working 2 h/day}{2}{5}
\end{ganttchart}
\caption{Either of the auditors working 8 or 4 h per day can be assigned to the depicted task. However, an auditor working only 2 h per day would not be able to finish the task in time.}
\label{fig:availability}
\end{figure}
}

\paragraph{Engagement availability}\label{p:engagementAvailability}
Each engagement has a calendar of days when auditors can work on the engagement. This allowed us to model the following constraints:
\e{Availability of the client data.}
\e{Due dates.}
\e{Client's company-wide holidays.}
\imp{An edge can be created only on a day when engagement $e$ is available. The edge can be created only if it is hypothetically possible for auditor $a$ to finish the task in time.}

\paragraph{Familiarity}\label{p:familiarity}
If possible, we prefer to assign auditors who have past experience with the engagement.
\imp{We solve this by adjusting the cost matrix $c$.}

\paragraph{Travel cost}
Auditor travel costs can have an important impact on the audit costs \cite{Chan1998}, particularly for audit companies with multiple offices. We approximate the travel cost by the flying distance between the auditor's office and the client in Euclidean projection and add it to the cost of the edge. 
Furthermore, not all auditors are willing to travel long distances.
\e{The auditor does not have a driving license.}
\imp{When an auditor is unwilling to travel, let us say, more than 10 km, we remove all the edges from the auditor to the clients that are more than 10 km away.

The travel cost between the office and the client is modeled by adjusting the cost matrix $c$.
}

\paragraph{Level substitutions}\label{p:levelSubstitutions}
 If necessary, an auditor at one level can be substituted by an auditor at a different level.
\e{If no junior analyst is available, assign a senior analyst on the junior analyst task}
\imp{We solve it by adjusting the cost matrix $c$.}

\paragraph{Updates}
The list of engagements and available auditors were not fixed. Hence, one of the requirements was the ability to change the engagement or the auditor lists without many changes to the already generated schedule. 
\imp{We approximate it by adjusting the cost matrix $c$ to include a penalty if an auditor will suddenly work on a new task. Conversely, we include a reward in the cost matrix if the auditor will work on the same task as before. Note that this formulation gives the optimizer a leeway to shift the tasks over time. 
}

\paragraph{Hard task preference}
The firm may require certain auditors to participate in certain tasks. While this is a trivial constraint, it is a frequent one. 
\e{The engagement manager must be assigned to his/her engagement.}
\e{Conflict of interests.}
\imp{
To enforce the assignment of auditor $a$ to task $t$, we set the following:
\begin{equation}
	\sum_{d \in D} x_{as} \ge 1 \text{ for auditor $a$ and all sinks $s$ with task $t$}
\end{equation}
The forbidden combinations are simply excluded from the set of candidate edges.
}

\paragraph{Soft task preference}
The auditors may express preferences to certain tasks \cite{Salewski1994}.
\e{The auditor wants to specialize in bank auditing.}
\imp{We solve this by adjusting the cost matrix $c$.}

\paragraph{Staff scarcity}\label{p:scarcity}
As the customer base is growing, there is an insufficient number of auditors to work all of the required hours. 
\imp{We solve this through the introduction of \q{mock auditors} in the set of auditors, who if assigned will have to be hired. Of course, the recruitment and training of new auditors is a costly endeavor. Hence, we penalize the employment of mock auditors by extending the objective function to:
\begin{equation}
	\underset{xy}{\text{minimize}} \sum_{a \in A} \sum_{s \in S} c_{as} x_{as} + c_{\text{mock}} \sum_{a \in M} y_a ,
\end{equation}
where $M$ \added{with} $M \subset A$ is the set of mock auditors, $c_{\text{mock}}$ is the cost per mock auditor, and $y_a$ is a binary slack variable for each mock auditor. We use following constraints to force slack variables $y_a$ to take value 1 when mock auditor $a$ is assigned to at least one sink $s$:
\begin{equation}
	y_a \geq x_{as} \text{ for each auditor $a \in M$ and sink $s$}.
\end{equation}
}

\paragraph{Uncertainty}
As we plan further into the future, the chance of unforeseen disturbances will increase. To minimize the expected count of schedule changes caused by unforeseen disturbances, we favor scheduling tasks earlier rather than later, given the choice.
\e{Some auditor leaves the firm.}
\imp{We implement this by adjusting the cost matrix $c$ with a hyperbolic discounting of the reward \cite{Sozou1998}:
\begin{equation}
	c_{\text{reward}}(d) = c_{\text{reward}} \frac{1}{1 + k d}, 
\end{equation}
where $c_{\text{reward}}$ is a reward for starting a task on the first day in the schedule, $d$ is the index of the day, and $k$ is a parameter governing the degree of discounting (we used $k = 0.01$ and constant $c_{\text{reward}}$).
}

\paragraph{Warm up}
It takes time to get ready to work on an engagement \cite{Salewski1996,Kadri2018}. This overhead is minimized by minimizing the number of distinct auditors at the engagement by assigning the same auditor to as many tasks as possible. To ensure that this criterion does not go against level specialization, this criterion uses a substantially smaller penalty than the level substitution penalty.
\imp{
We extend the objective function to the following:
\begin{equation}
	\underset{xyz}{\text{minimize}} \sum_{a \in A} \sum_{s \in S} c_{as}x_{as} + c_{\text{mock}} \sum_{a \in M} y_a + c_{\text{warm up}} \sum_{a \in A} \sum_{e \in E} z_{ae},
\end{equation}
where $c_{\text{warm up}}$ is the cost of the \q{warming up} of one auditor to one engagement, and $z_{ae}$ is a binary slack variable for each auditor and engagement. We use following constraints to force slack variables $z_{ae}$ to take value 1 when auditor $a$ is assigned to engagement $e$ at least once:
\begin{equation}
	z_{ae} \geq x_{as} \text{ for auditor $a$ and all sinks $s$ with engagement $e$}.
\end{equation}
}

\subsection{Rejected constraints}
We also present a list of constraints, which were rejected by our auditors, but which, nevertheless, might be of interest to the reader.

\paragraph{Precedence} Task $j$ cannot start before task $i$ is completed \cite{Brucker2001}. 
\e{The audit must be concluded before the documents get translated.}
\just{Our auditors rejected this constraint because it is already considered in the \hyperref[p:engagementAvailability]{engagement availability} windows.}
  
 \paragraph{Lag} Minimum and maximum time-lags between two tasks are given \cite{Salewski1994}.
\e{By law, the period between two tasks cannot exceed 3 months.}
\just{Our auditors rejected this constraint because it is already considered in the \hyperref[p:engagementAvailability]{engagement availability} windows.}

\paragraph{Parallelity} Tasks may be forced to be processed in parallel \cite{Brucker2001}.
\e{An analyst and the analyst's manager must work on the engagement on the same days.}
\just{Our auditors did not consider this to be an actual problem to solve.}

\paragraph{Workload} The total workload of an auditor is lower bounded \cite{Brucker2001}.
\e{An auditor must have at least 30\% utilization.}
\just{The firm does not suffer from the underutilization of the auditors but rather from the opposite problem. Consequently,  \hyperref[p:auditorAvailability]{auditor availability} constraints suffice.}

\section{Implementation}\label{sec:implementation}
The code is written in Python 3.7 and is divided into 3 parts:
\begin{description}
  \item[Load] Data loading.
  \item[Optimization] Problem formulation and optimization.
  \item[Export] Solution extraction.
\end{description}

\begin{figure*}[htb]
	\centering
	\includegraphics[width=\textwidth]{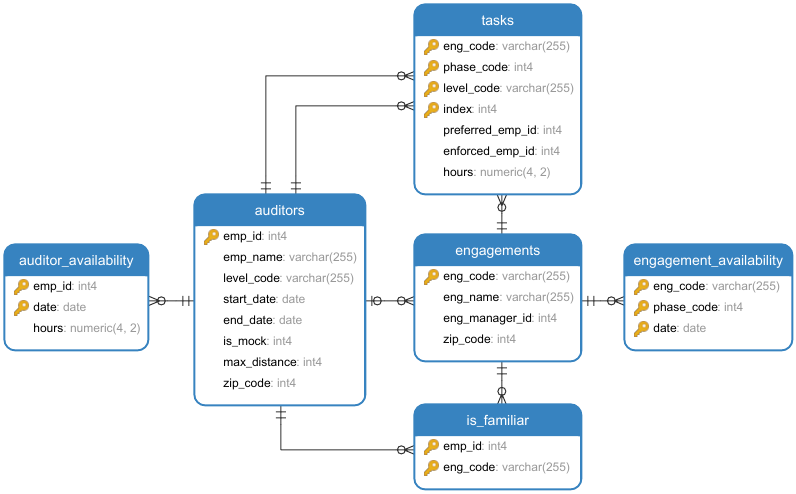}
	\caption{Entity-relationship diagram of the input data in IE notation.}
	\label{fig:er}
\end{figure*}

The input data are in an SQL database. The advantage of using the SQL database as a data storage is that it allows us to define integrity constraints (e.g., not null, unique, foreign key) declaratively and that essentially each data engineer is familiar enough with SQL to deal with eventual integrity violations. The entity-relationship diagram of the input data is in \cref{fig:er}.






\added{The input data, which must be provided by the auditors for a new schedule, include a set of \texttt{auditors}, a set of \texttt{engagements}, availability of the auditors and engagements (\texttt{auditor\_availability} and \texttt{engagement\_availability} respectively), a set of \texttt{tasks} that the auditors should do, and a set of the past auditor assignments to the engagements as a proxy of the auditors' familiarity with the engagements (\texttt{is\_familiar}).

Inputs not listed in \cref{fig:er} include parameters, like the travel cost in EUR for 1 km, which are used to calculate the cost matrix $C$ from the input data, and a potential past solution $X$ that should be updated to conform to the new requirements.}

For optimization, we experimentally compared open-source OR-Tools v7.7.7810 (with CBC v2.10.5 solver in the backend) to commercial Gurobi v9.0.3.

\section{Results} \label{sec:results}
The implementation was evaluated on real data with 71 auditors, 47 engagements, 3 phases, 10 levels, 6 indices, and 365 days. The runtime measurements were obtained on a laptop with an Intel Core i5 processor. 

An empirical comparison of the computer-generated schedule to the handmade schedule by an experienced planner based on the criteria listed in \cref{sec:requirements} is presented in \cref{tab:comparison}.

\begin{table}[htb]
\centering
\begin{tabular}{lrr}
	\toprule
	Criterium & Handmade & Computer-generated \\
	\midrule
	\# \nameref{p:auditorAvailability} violations & 4 & 0 \\
	\# \nameref{p:engagementAvailability} violations & 46 & 0 \\
	\# \nameref{p:familiarity} violations & 157 & 153 \\
	\# \nameref{p:levelSubstitutions} & 0 & 41 \\
	\# \hyperref[p:scarcity]{Auditors to recruit} & 3 & 2 \\
	Time to collect data & 38 hours & 38 hours \\
	Time to schedule & 24 hours & \num{2654} seconds \\
	Time to update schedule & 8 hours & \num{236} seconds \\
	\bottomrule
\end{tabular}
\caption{Empirical comparison of a handmade schedule to the optimized schedule. When a handmade schedule takes 24 hours to generate, it actually means 3 man-days.}
\label{tab:comparison}
\end{table}

The computer-generated schedule is better or equal to the handmade schedule in all regards except the count of level substitutions. To be fair to the experienced planner and to further highlight the difficulty of finding a feasible schedule given the constraints, if we did not permit level substitutions in the computer-generated schedule, a feasible solution (a solution that does not violate any hard constraint) would require at least five new auditors owing to the clustering of engagements over time.

The comparison also includes the time to collect the engagement availability from the clients and the task set from the engagement managers (based on the survey conducted, each manager spends approximately 2 h on this task). We do not use any additional information beyond what has already been collected.
The reported runtime of \num{2654} s includes input data validation \& preprocessing, model creation \& optimization, and schedule export. This runtime was a\deleted{r}chived using the Gurobi solver. When we used OR-Tools, the runtime increased to \num{55793} s (21-times more than with Gurobi).
The runtime of the schedule update (e.g., adjustment of the availability of an auditor or an engagement) is dominated by the model creation as the solver is warm-started from the previous solution.

\section{Discussion} \label{sec:discussion}
Since this article is not only about solving the audit scheduling problem exactly but also about solving the problem quickly, we have to discuss implementation details.

\paragraph{Hard vs. soft constraints}
In comparison to the reference ILP formulation of audit scheduling by Balachandran and Zoltner \cite{Zoltners1989}, we exclude the illegal assignments from the problem formulation, whereas Balachandran and Zoltner leave them in the problem formulation and only adjust the corresponding cost in $c_{as}$ to a large constant. The advantage of exclusion is a smaller memory footprint and faster optimization (observe \cref{tab:soft}).
\e{
When we conduct an exhaustive cross-join between 71 auditors, 47 engagements, 3 phases, 10 levels, 6 indices, and 365 days, we obtain \num{219240900} edges. However, in real data of the same size, only \num{245455} are needed (an almost \num{1000}-fold reduction).
}

\begin{table}[htb]
\centering
\begin{tabular}{lrr}
\toprule
Implementation & Peak memory {[}MB{]} & Runtime {[}s{]} \\
\midrule
Hard constraint & \num{61} & \num{2654} \\
Soft constraint & \num{32311} & \num{132649} \\
\bottomrule
\end{tabular}
\caption{The effect of hard vs. soft constraints.}
\label{tab:soft}
\end{table}

\paragraph{Early pruning}
While modern ILP implementations have efficient presolvers that can quickly prune away unused or fixed variables, we found early pruning during the problem formulation highly beneficial, because it leads to smaller matrices and that speeds up all the following data transformations (observe \cref{tab:pruning}). Examples of early pruning include exclusion of edges where an auditor cannot finish the workload in time, or exclusion of days when auditors do not work, like national holidays, as presented in \cref{p:auditorAvailability}.

\begin{table}[htb]
\centering
\begin{tabular}{lrr}
\toprule
Implementation & Peak memory {[}MB{]} & Runtime {[}s{]} \\
\midrule
With early pruning & \num{61} & \num{2654} \\
Without early pruning & \num{630} & \num{2791} \\
\bottomrule
\end{tabular}
\caption{The effect of early pruning.}
\label{tab:pruning}
\end{table}

\paragraph{Problem representation}
A tempting problem representation of the scheduling problem is to use a binary vector $x_{aed}$:

\begin{equation}
	x_{aed} =
	\begin{cases}
		1 & \text{if auditor $a$ works on engagement $e$ on day $d$} \\
	    0 & \text{otherwise.}
	\end{cases}
\end{equation}

However, similarly to Chan, Dodin et al. \cite{Chan1986,Dodin1991,Dodin1997,Chan1998}, we opted to use a binary vector \replaced{$x_{aed}$}{$x_{a,e,d}$}:
\begin{equation}
	x_{aed} =
	\begin{cases}
		1 & \text{if auditor $a$ started to work on engagement $e$ on day $d$,} \\
	    0 & \text{otherwise,}
	\end{cases}
\end{equation}
because it gives us \added{the} continuity constraint \added{of} \cref{p:continuity} for free without the need to refine the optimization criterium to also minimize the length of the assignments (observe the impact in \cref{tab:representation}).

\begin{table}[htb]
\centering
\begin{tabular}{lrr}
\toprule
Implementation & Peak memory {[}MB{]} & Runtime {[}s{]} \\
\midrule
\replaced{Auditor started to work}{Auditor works} & \num{61} & \num{2654} \\
\replaced{Auditor works}{Auditor started to work} & \num{90} & \num{16493} \\
\bottomrule
\end{tabular}
\caption{The effect of changing what $x$ means.}
\label{tab:representation}
\end{table}

\paragraph{Mock auditors}
The advantage of using a large pool of mock auditors is that the solver can find an initial solution quickly. The disadvantage is that it then takes longer to converge to the optimal solution and to prove that the solution is indeed optimal (observe \cref{tab:mock}). 

\begin{table}[htb]
\centering
\begin{tabular}{lrr}
\toprule
Implementation & Peak memory {[}MB{]} & Runtime {[}s{]} \\
\midrule
2 mock auditors & 61 & \num{2713} \\
3 mock auditors & 61 & \num{2654} \\
10 mock auditors & 62 & \num{2701} \\
20 mock auditors & 62 & \num{2816} \\
\bottomrule
\end{tabular}
\caption{The effect of the number of mock auditors.}
\label{tab:mock}
\end{table}

We recommend using just a slightly higher number of mock auditors than is the expected need to minimize the runtime.

\paragraph{Divide and conquer}
After the initial run on a subset of 71 auditors and 47 engagements, described in \cref{sec:results}, the auditing company decided to generate the schedule for all 299 auditors and 271 engagements. This led to an increase in memory requirements and runtime. But after dividing auditors and engagements into 4 \replaced{disjoint}{disjunct} subsets, particularly memory requirements decreased significantly (observe \cref{tab:divide}). Since the division was based on distinct types of audit (banking, non-banking) and geographically distinct areas (east, west), it did not lead, at least in this case, to a degradation of the obtained solution.

\begin{table}[htb]
\centering
\begin{tabular}{lrr}
\toprule
Implementation & Peak memory {[}MB{]} & Runtime {[}s{]} \\
\midrule
Whole at once & \num{44316} & \num{23868} \\
4 subsets one after another & \num{92} & \num{8292} \\
\bottomrule
\end{tabular}
\caption{The effect of dividing the problem into subproblems.}
\label{tab:divide}
\end{table}

\section{Conclusion} \label{sec:conclusion}

We implemented a schedule optimizer for an auditing firm and published the code at \href{https://github.com/janmotl/audit-scheduling}{github.com/janmotl/audit-scheduling}. We formulated the scheduling problem as a multi-commodity network flow problem. Our contribution is in adapting a multi-commodity network flow formulation from \textit{shift-based} scheduling with fixed starting and end times, to \textit{task-based} scheduling, where the tasks can be performed at any time in a temporal window. 

According to the auditors, the introduction of some automated planners was necessary because the handmade scheduling started to take an unacceptably long time. Our planner implementation reduced the time to schedule from 3 days to 45 min. As a byproduct, the computer-generated schedule is of higher quality than the handmade schedule. The most appreciated improvement was the reduction in the number of auditors to recruit from three to two, because the wage makes a significant part of all their expenses.

In \deleted{a} future \added{work}, we plan to extend the planner to consider the flu season when some of the auditors are expected to be unavailable.  

\bibliography{library}{}
\bibliographystyle{plain}

\end{document}